\documentclass[11pt,a4paper]{article}

\usepackage{graphics,epsfig,theorem,latexsym,amssymb,amsmath,amsfonts}
\usepackage{times,epic,eepic}

\usepackage{pifont}

\usepackage{cite}
\usepackage{url}

\usepackage{color}

\usepackage{newlfont}

 \newenvironment{keywords}{ \noindent {\small\bf Key Words}:}{ }

\def\G1{\hbox{$\displaystyle{\mbox{\ding{172}}}$}}

\def\bd{\begin{description}}
\def\ed{\end{description}}

\def\beq{\begin{equation}}
\def\eeq{\end{equation}}
\def\bea{\begin{eqnarray}}
\def\eea{\end{eqnarray}}
\def\beas{\begin{eqnarray*}}
\def\eeas{\end{eqnarray*}}

\newtheorem{corollary}{Corollary}
\newtheorem{theorem}{Theorem}

\begin{document}

\title{  Solving ordinary differential equations\\ on the Infinity Computer\\
by working with infinitesimals numerically\thanks{This study was
supported by the Ministry of
   Education and Science of Russian Federation, project
   14.B37.21.0878. The author thanks anonymous reviewers for their useful suggestions.}}

\newcommand{\nms}{\normalsize}
\author{\\ { \vspace*{-5mm}  \bf Yaroslav D. Sergeyev\footnote{Yaroslav D.
Sergeyev, Ph.D., D.Sc., is   Distinguished
 Professor  at the University of Calabria, Rende, Italy.
 He is also Full Professor (part-time contract) at the N.I.~Lobatchevsky State University,
  Nizhni Novgorod, Russia and   Affiliated  Researcher at the Institute of High Performance Computing
   and Networking of the National Research Council of Italy.
 }\,\,
    }\\[-4pt]
    \\
      \nms Dipartimento di Elettronica, Informatica e Sistemistica,\\[-4pt]
       \nms   Universit\`a della Calabria,\\[-4pt]
       \nms 87030 Rende (CS)  -- Italy\\ \\[-4pt]
        \nms tel./fax: +39 0984 494855\\[-4pt]
       \nms http://wwwinfo.deis.unical.it/$\sim$yaro\\[-4pt]
         \nms {\tt  yaro@si.deis.unical.it }
}

\date{}

\maketitle

\vspace{-8mm}

 \begin{abstract}
 There exists a huge number
of numerical methods that   iteratively construct approximations to
the solution $y(x)$ of an ordinary  differential equation (ODE)
$y'(x)=f(x,y)$ starting from an initial value $y_0=y(x_0)$ and using
a finite approximation step $h$ that influences the accuracy of the
obtained approximation. In this paper, a new framework for solving
ODEs is presented for a new kind of a computer -- the Infinity
Computer
 (it has been patented and    its working
prototype exists). The new computer is  able to work numerically
with finite, infinite, and infinitesimal numbers giving so the
possibility to use different infinitesimals numerically and, in
particular, to take advantage of infinitesimal   values of~$h$. To
show the potential of the new framework a number of results is
established. It is proved that the Infinity Computer is able to
calculate derivatives of the solution $y(x)$ and to reconstruct its
Taylor expansion of a desired order numerically  without   finding
the respective derivatives analytically (or symbolically) by the
successive derivation of the ODE   as it is usually done when the
Taylor method is applied. Methods using approximations of
derivatives obtained thanks to infinitesimals are discussed  and a
technique for an automatic control of rounding errors   is
introduced. Numerical examples are given.

 \end{abstract}

\vspace{3mm}



\vspace{5mm}

\section{Introduction}
\label{s_m1}

The number of  applications in physics, mechanics, and engineering
where it is necessary to solve numerically ordinary differential
equations (ODEs) with a given initial value is really enormous.
Since many ordinary differential equations   cannot be solved
analytically, people use numerical algorithms\footnote{ There exist
also symbolic techniques but they are not considered in this paper
dedicated to numerical computations.} for finding approximate
solutions (see \cite{Brugnano,Butcher,Henrici,Quarteroni}). In this
paper, we want to approximate the solution $y(x),\,\, x \in [a,b],$
of the initial value problem (also called the Cauchy problem) for a
differential equation
 \beq
 y'(x)=f(x,y), \hspace{1cm} y(x_0)=y_0,
 \hspace{1cm} x_0 = a, \label{ode1}
 \eeq
where $a$ and $b$ are finite numbers and $y(x_0)=y_0$ is called the
initial condition. We suppose that $f(x,y)$ is given by a computer
procedure.  Since very often in scientific and technical
applications it can happen that the person who wants to solve
(\ref{ode1})   is not the person who has written the code for
$f(x,y)$, we suppose that the person solving (\ref{ode1}) does not
know the structure of  $f(x,y)$, i.e., it is a black box for
him/her.

In the literature, there exist numerous numerical algorithms
constructing a sequence $y_1, y_2, y_3, \ldots$ approximating the
exact values $y(x_1), y(x_2), y(x_3), \ldots$ that the solution
$y(x)$ assumes at points $x_1, x_2, x_3, \ldots$ (see
\cite{Butcher,Hairer,Lambert}). The explicit Euler algorithm is the
simplest among explicit methods for the numerical integration of
ODEs. It uses the first two terms of the Taylor expansion of $y(x)$
constructing so the linear approximation around the point
$(x_0,y(x_0))$. The $(n+1)$th step of the Euler algorithm describes
how to move  from the point $x_n$ to $x_{n+1} = x_{n} + h, n>0,$ and
is executed as follows
 \beq
y_{n+1} = y_n + h f(x_n,y_n). \label{ode2}
 \eeq

Traditional computers work with finite values of $h$ introducing so
errors at each step of the algorithm.  In order to obtain more
accurate approximations it is necessary to decrease the step $h$
increasing so the number of steps of the method (the computations
become more expensive). In any case, $h$ always remains finite and
its minimal acceptable value is determined by technical
characteristics of each concrete computer the method is implemented
on. Obviously, the same effects hold for more sophisticated methods,
as well (see \cite{Butcher,Hairer,Henrici,Lambert}). Another
approach to solve (\ref{ode1}) on a traditional computer is the use
of an automatic differentiation software executing pre-processing of
(\ref{ode1}) (see \cite{Griewank_Corliss} and references given
therein).

In this paper, we introduce a new numerical framework for solving
ODEs related to the usage of a new kind of computer -- the Infinity
Computer (see \cite{Sergeyev_patent,www,informatica}). It  is able
to work \textit{numerically} with finite, infinite, and
infinitesimal quantities. The Infinity Computer is based on   an
applied point of view (see \cite{Sergeyev,informatica,Lagrange}) on
infinite and infinitesimal numbers.   In order to see the place of
the new approach in the historical panorama of ideas dealing with
infinite and infinitesimal, see
\cite{Lolli,MM_bijection,Dif_Calculus,first,Sergeyev_Garro}. The new
methodology has been successfully applied for studying numerical
differentiation and optimization (see
\cite{DeLeone,Korea,Num_dif,Zilinskas}), fractals (see
\cite{chaos,Menger,Biology,DeBartolo}), percolation (see
\cite{Iudin,DeBartolo}), Euclidean and hyperbolic geometry (see
\cite{Margenstern,Rosinger2}), the first Hilbert problem and Turing
machines (see \cite{first,Sergeyev_Garro,Sergeyev_Garro_2}),
cellular automata (see \cite{DAlotto}), infinite series (see
\cite{Dif_Calculus,Riemann,Zhigljavsky}), functions and their
derivatives that can assume infinite and infinitesimal values (see
\cite{Dif_Calculus}), etc.

With respect to the initial value problem (\ref{ode1}), the
possibility to work numerically with infinitesimals allows us to use
\textit{numerical infinitesimal values of $h$}. It is proved that
under reasonable conditions the Infinity Computer is able to
calculate \textit{exact} values of the derivatives of $y(x)$ and to
reconstruct its Taylor expansion with a desired accuracy by using
 infinitesimal values of $h$  without   finding
the respective derivatives analytically (or symbolically) by the
successive derivation of (\ref{ode1}) as it is usually done when the
Taylor method is applied.

The rest of the paper is organized as follows. Section~2 briefly
presents the new computational methodology. Section~3 introduces the
main theoretical results and describes how derivatives of $y(x)$ can
be calculated numerically on the Infinity Computer.  Section~4
introduces a variety of examples of the usage of infinitesimals for
ODEs numerical solving. First, it presents two simple iterative
methods. Then, it describes a technique  that can be used to obtain
approximations of derivatives of the solution $y(x)$ at the point
$x_{n+1}$ using infinitesimals and the information obtained at the
point $x_{n}$. Finally, a technique for an automatic control of
rounding errors that can occur  during evaluation of  $f(x,y)$  is
introduced. Through the paper, theoretical results are illustrated
by numerical examples.

\section{A fast tour to the new computational methodology}

 Numerous trials
have been done during the centuries in order to evolve existing
numeral systems\footnote{ We are reminded that a \textit{numeral} is
a symbol or group of symbols that represents a \textit{number}. The
difference between numerals and numbers is the same as the
difference between words and the things they refer to. A
\textit{number} is a concept that a \textit{numeral} expresses. The
same number can be represented by different numerals. For example,
the symbols `7', `seven', and `VII' are different numerals, but they
all represent the same number.} in such a way that infinite and
infinitesimal numbers could be included in them (see
\cite{Benci,Cantor,Conway,Leibniz,Levi-Civita,Newton,Robinson,Wallis}).
Particularly, in the early history of the calculus, arguments
involving infinitesimals played a pivotal role in the derivation
developed by Leibniz and Newton (see \cite{Leibniz,Newton}). The
notion of an infinitesimal, however, lacked a precise mathematical
definition and in order  to provide a more rigorous   foundation for
the calculus, infinitesimals were gradually replaced by the
d'Alembert-Cauchy  concept of a limit.

Since new numeral systems appear very rarely, in each concrete
historical period their importance for Mathematics is very often
underestimated (especially by pure mathematicians). In order to
illustrate their importance, let us remind the Roman numeral system
that does not allow one to express zero and negative numbers. In
this system, the expression III-X is an indeterminate form. As a
result, before appearing the positional numeral system and inventing
zero (by the way, the second event was several hundred years later
with respect to the first one) mathematicians were not able to
create theorems involving zero and negative numbers and to execute
computations with them.

There exist    numeral systems that are even weaker than the Roman
one. They seriously limit  their users in  executing computations.
Let us recall a study published recently in \textit{Science} (see
\cite{Gordon}) that describes a primitive tribe -- Pirah\~{a} --
living in Amazonia. These people use a very simple numeral system
for counting: one, two, many. For Pirah\~{a}, all quantities larger
than two are just `many' and such operations as 2+2 and 2+1 give the
same result, i.e., `many'. Using their weak numeral system
Pirah\~{a} are not able to see, for instance, numbers 3, 4, 5, and
6, to execute arithmetical operations with them, and, in general, to
say anything about these numbers because in their language there are
neither words nor concepts for that.

In the context of the present paper, it is very important that  the
weakness of Pirah\~{a}'s numeral system leads them to such results
as
 \beq
 \mbox{`many'}+ 1= \mbox{`many'},   \hspace{1cm}
\mbox{`many'} + 2 = \mbox{`many'},
  \label{piraha1}
       \eeq
which are very familiar to us  in the context of views on infinity
used in the traditional calculus
 \beq
  \infty + 1= \infty,
\hspace{1cm}    \infty + 2 = \infty.
 \label{piraha2}
       \eeq
The arithmetic of Pirah\~{a} involving the numeral `many' has also a
clear similarity with the arithmetic proposed by Cantor for his
Alephs\footnote{This similarity becomes even more pronounced  if one
considers another Amazonian tribe -- Munduruk\'u (see \cite{Pica})
-- who fail in exact arithmetic with numbers larger than 5 but are
able to compare and add large approximate numbers that are far
beyond their naming range. Particularly, they use the words `some,
not many' and `many, really many' to distinguish two types of large
numbers using the rules that are very similar to ones used by Cantor
to operate with $\aleph_0$ and $\aleph_1$, respectively.}:
 \beq
\aleph_0 + 1= \aleph_0,    \hspace{1cm} \aleph_0 + 2= \aleph_0,
\hspace{1cm}\aleph_1+ 1= \aleph_1,    \hspace{1cm}   \aleph_1 + 2 =
\aleph_1.
  \label{piraha3}
       \eeq

Thus, the modern mathematical numeral systems allow us to
distinguish a larger quantity of finite numbers with respect to
Pirah\~{a} but give results that are   similar to those of
Pirah\~{a} when we speak about infinite numbers. This observation
leads us to the following idea: \textit{Probably our difficulties in
working with infinity is not connected to the nature of infinity
itself but is a result of inadequate numeral systems that we use to
work with infinity, more precisely, to express infinite numbers.}

Let us compare the usage of numeral systems in Mathematics
emphasizing differences that hold when one works, on the one hand,
with finite quantities and, on the other hand, with infinities and
infinitesimals. In our every day activities with finite numbers the
\emph{same} finite numerals are used for \emph{different} purposes
(e.g., the same numeral 4 can be used to express   the number of
elements of a set and to indicate the position of an element in a
finite sequence). When we face the necessity to work with infinities
or infinitesimals, the situation changes drastically. In fact, in
this case \emph{different} symbols are used to work with infinities
and infinitesimals in \emph{different} situations:
\begin{itemize}
\item
 $\infty$ in standard Analysis;
 \item
$\omega$ for working with ordinals;
\item
 $\aleph_0, \aleph_1, ...$
for dealing with cardinalities;
\item
 non-standard numbers using a
generic infinitesimal $h$ in non-standard Analysis, etc.
\end{itemize}

In particular,   since the mainstream of the traditional Mathematics
very often does not pay any attention to the distinction between
numbers and numerals (in this occasion it is necessary to recall
constructivists who studied this issue), many theories dealing with
infinite and infinitesimal quantities have a symbolic (not
numerical) character. For instance, many versions of the
non-standard Analysis are symbolic, since they have no numeral
systems to express their numbers by a finite number of symbols (the
finiteness of the number of symbols  is necessary for organizing
numerical computations). Namely, if we consider a finite $n$ than it
can be taken $n=5$, or $n=103$ or any other numeral used to express
finite quantities and consisting of a finite number of symbols. In
contrast, if we consider a non-standard infinite $m$ then it is not
clear which numerals can be used to assign a concrete value to $m$.

 Analogously, in non-standard
Analysis, if we consider an infinitesimal $h$ then it is not clear
which numerals consisting of a finite number of symbols can be used
to assign a value to $h$ and to write $h=...$ In fact, very often in
non-standard Analysis texts, a \textit{generic} infinitesimal $h$ is
used and it is considered as a symbol, i.e., only symbolic
computations can be done with it. Approaches of this kind leave
unclear such issues, e.g., whether the infinite $1/h$ is integer or
not or whether $1/h$ is the number of elements of an infinite set.
Another problem is related to comparison of values. When we work
with finite quantities then we can compare $x$ and $y$ if they
assume numerical values, e.g., $x=4$ and $y=6$ then, by using rules
of the numeral system the symbols 4 and 6 belong to, we can compute
that $y>x$. If one wishes to consider two infinitesimals $h_1$ and
$h_2$ then it is not clear how to compare them because numeral
systems that can express infinitesimals are not provided by
non-standard Analysis techniques.

 The approach developed in
\cite{Sergeyev,informatica,Lagrange} proposes a numeral system that
uses \textit{the same numerals}  for several different purposes for
dealing with infinities and infinitesimals: in Analysis for working
with functions that can assume different infinite, finite, and
infinitesimal values (functions can also have derivatives assuming
different infinite or infinitesimal values); for measuring infinite
sets; for indicating positions of elements in ordered infinite
sequences; in probability theory, etc. It is important to emphasize
that the new numeral system   avoids situations of the type
(\ref{piraha1})--(\ref{piraha3}) providing results ensuring that  if
$a$ is a numeral written in this system then for any $a$ (i.e., $a$
can be finite, infinite, or infinitesimal)   it follows $a+1>a$.

The new numeral system works as follows. A new infinite unit of
measure expressed by the numeral \G1 called \textit{grossone} is
introduced as the number of elements of the set, $\mathbb{N}$, of
natural numbers. Concurrently with the introduction of grossone in
the mathematical language all other symbols (like $\infty$, Cantor's
$\omega$, $\aleph_0, \aleph_1, ...$,  etc.) traditionally used to
deal  with infinities and infinitesimals are excluded from the
language because grossone and other numbers constructed with its
help not only can be used instead of all of them but can be used
with a higher accuracy. Grossone is introduced by describing its
properties postulated by the Infinite Unit Axiom (see
\cite{informatica,Lagrange}) added to axioms for real numbers
(similarly, in order to pass from the set, $\mathbb{N}$, of natural
numbers to the set, $\mathbb{Z}$, of integers a new element -- zero
expressed by the numeral 0 -- is introduced by describing its
properties).

 The new
numeral \G1 allows us to construct different numerals expressing
different infinite and infinitesimal numbers and to execute
computations with them. As a result, in Analysis, instead of the
usual symbol $\infty$ used in series and integration different
infinite and/or infinitesimal numerals can be used (see
\cite{Dif_Calculus,Riemann,Zhigljavsky}). Indeterminate forms are
not present and, for example, the following relations hold for
$\mbox{\ding{172}}$ and $\mbox{\ding{172}}^{-1}$ (that is
infinitesimal), as for any other (finite, infinite, or
infinitesimal) number expressible in the new numeral system
 \beq
 0 \cdot \mbox{\ding{172}} =
\mbox{\ding{172}} \cdot 0 = 0, \hspace{3mm}
\mbox{\ding{172}}-\mbox{\ding{172}}= 0,\hspace{3mm}
\frac{\mbox{\ding{172}}}{\mbox{\ding{172}}}=1, \hspace{3mm}
\mbox{\ding{172}}^0=1, \hspace{3mm} 1^{\mbox{\tiny{\ding{172}}}}=1,
\hspace{3mm} 0^{\mbox{\tiny{\ding{172}}}}=0,
 \label{3.2.1}
       \eeq
\[
 0 \cdot \mbox{\ding{172}}^{-1} =
\mbox{\ding{172}}^{-1} \cdot 0 = 0, \hspace{5mm}
\mbox{\ding{172}}^{-1} > 0, \hspace{5mm} \mbox{\ding{172}}^{-2} > 0,
\hspace{5mm} \mbox{\ding{172}}^{-1}-\mbox{\ding{172}}^{-1}= 0,
\]
\[
\frac{\mbox{\ding{172}}^{-1}}{\mbox{\ding{172}}^{-1}}=1,
\hspace{3mm}
\frac{\mbox{\ding{172}}^{-2}}{\mbox{\ding{172}}^{-2}}=1,
\hspace{3mm} (\mbox{\ding{172}}^{-1})^0=1, \hspace{5mm}
\mbox{\ding{172}} \cdot \mbox{\ding{172}}^{-1} =1, \hspace{5mm}
\mbox{\ding{172}} \cdot \mbox{\ding{172}}^{-2}
=\mbox{\ding{172}}^{-1}.
       \]

The new approach gives the possibility to develop a new Analysis
(see  \cite{Dif_Calculus}) where  functions assuming not only finite
values but also infinite and infinitesimal ones can be studied. For
all of them it becomes possible to introduce a new notion of
continuity that is closer to our modern physical knowledge.
Functions  assuming finite and infinite values can be differentiated
and integrated.

\textbf{Example 1.} \label{e_m1} The function $f(x)=x^2$ has the
first derivative $f'(x)= 2x$ and both $f(x)$ and $f'(x)$ can be
evaluated at infinite and infinitesimal $x$. Thus, for infinite
$x=\mbox{\ding{172}}$ we obtain infinite values
\[
f(\mbox{\ding{172}})= \mbox{\ding{172}}^{2}, \hspace{1cm}
f'(\mbox{\ding{172}})= 2\mbox{\ding{172}}
\]
and for infinitesimal $x=\mbox{\ding{172}}^{-1}$ we have
infinitesimal values
\[
f(\mbox{\ding{172}}^{-1})= \mbox{\ding{172}}^{-2}, \hspace{1cm}
f'(\mbox{\ding{172}}^{-1})= 2\mbox{\ding{172}}^{-1}.
\]
If $x=5\mbox{\ding{172}}-10\mbox{\ding{172}}^{-1}$ then we have
\[
f(\mbox{\ding{172}}^{-1})=
(5\mbox{\ding{172}}-10\mbox{\ding{172}}^{-1})^2=
25\mbox{\ding{172}}^{2}-100 +100\mbox{\ding{172}}^{-2},
 \]
 \[
f'(\mbox{\ding{172}}^{-1})=
10\mbox{\ding{172}}-20\mbox{\ding{172}}^{-1}.
\]
 We can also work with functions defined by formulae
including infinite and infinitesimal numbers.  For example, the
function $f(x)=\frac{1}{\mbox{\ding{172}}}x^2+\mbox{\ding{172}}x$
has a quadratic term infinitesimal and the linear one infinite. It
has the first derivative $f'(x)= \frac{2}{\mbox{\ding{172}}}x
+\mbox{\ding{172}}$. For infinite $x=3\mbox{\ding{172}}$ we obtain
infinite values
\[
f(\mbox{\ding{172}})= 3\mbox{\ding{172}}^{2} + 9\mbox{\ding{172}},
\hspace{1cm}f'(\mbox{\ding{172}})= \mbox{\ding{172}}+6
\]
and for infinitesimal $x=\mbox{\ding{172}}^{-1}$ we have
 \[
\hspace{20mm} f(\mbox{\ding{172}}^{-1})= 1+\mbox{\ding{172}}^{-3},
\hspace{1cm}
 f'(\mbox{\ding{172}}^{-1})= \mbox{\ding{172}}+
2\mbox{\ding{172}}^{-2}. \hspace{20mm} \Box
 \]

 By using the new numeral system it becomes   possible
to measure certain infinite sets and   to see, e.g., that the sets
of even and odd numbers have $\G1/2$ elements each. The set,
$\mathbb{Z}$, of integers has $2\G1+1$ elements (\G1 positive
elements, \G1 negative elements, and zero). Within the countable
sets and sets having cardinality of the continuum (see
\cite{Lolli,first,Lagrange}) it becomes possible to distinguish
infinite sets having different number of elements expressible in
the numeral system using grossone and to see that, for instance,
 \[
 \frac{\mbox{\ding{172}}}{2} < \mbox{\ding{172}}-1 < \mbox{\ding{172}} < \mbox{\ding{172}}+1 < 2\mbox{\ding{172}}+1 <
 2\mbox{\ding{172}}^2-1 < 2\mbox{\ding{172}}^2    <
 2\mbox{\ding{172}}^2+1  <
\]
 \[
2\mbox{\ding{172}}^2+2  <  2^{\mbox{\ding{172}}}-1 <
2^{\mbox{\ding{172}}} < 2^{\mbox{\ding{172}}}+1 <
10^{\mbox{\ding{172}}} <
  \mbox{\ding{172}}^{\mbox{\ding{172}}}-1 <
  \mbox{\ding{172}}^{\mbox{\ding{172}}} <
  \mbox{\ding{172}}^{\mbox{\ding{172}}}+1.
  \]

The Infinity Computer used in this paper for solving the problem
(\ref{ode1}) works with numbers having finite, infinite, and
infinitesimal parts. To represent them in the computer memory
records similar to traditional positional numeral systems can be
used (see \cite{Sergeyev_patent,informatica}). To construct a number
$C$ in the new numeral positional system\footnote{ At the first
glance the numerals (\ref{3.12}) can remind  numbers from the
Levi-Civita field (see \cite{Levi-Civita}) that is a very
interesting and important precedent of   algebraic manipulations
with infinities and infinitesimals. However, the two mathematical
objects have several crucial differences. They have been introduced
for different purposes by using two   mathematical languages having
different accuracies and on the basis of different methodological
foundations. In fact, Levi-Civita does not discuss the distinction
between numbers and numerals. His numbers have neither cardinal nor
ordinal properties; they are build using a generic infinitesimal and
only its rational powers are allowed; he uses symbol $\infty$ in his
construction; there is no any numeral system that would allow one to
assign numerical values to these numbers; it is not explained how it
would be possible to pass from d a generic infinitesimal~$h$ to a
concrete one  (see also the discussion above on the distinction
between numbers and numerals). In no way the said above should  be
considered as a criticism  with respect to results of Levi-Civita.
The above discussion has been introduced in this text  just to
underline that we are in front of two different mathematical tools
that   should be used in different mathematical contexts. } with
base \ding{172}, we subdivide $C$ into groups corresponding to
powers of \ding{172}:
 \beq
  C = c_{p_{m}}
\mbox{\ding{172}}^{p_{m}} +  \ldots + c_{p_{1}}
\mbox{\ding{172}}^{p_{1}} +c_{p_{0}} \mbox{\ding{172}}^{p_{0}} +
c_{p_{-1}} \mbox{\ding{172}}^{p_{-1}}   + \ldots   + c_{p_{-k}}
 \mbox{\ding{172}}^{p_{-k}}.
\label{3.12}
       \eeq
 Then, the record
 \beq
  C = c_{p_{m}}
\mbox{\ding{172}}^{p_{m}}    \ldots   c_{p_{1}}
\mbox{\ding{172}}^{p_{1}} c_{p_{0}} \mbox{\ding{172}}^{p_{0}}
c_{p_{-1}} \mbox{\ding{172}}^{p_{-1}}     \ldots c_{p_{-k}}
 \mbox{\ding{172}}^{p_{-k}}
 \label{3.13}
       \eeq
represents  the number $C$, where all numerals $c_i\neq0$, they
belong to a traditional numeral system and are called
\textit{grossdigits}. They express finite positive or negative
numbers and show how many corresponding units
$\mbox{\ding{172}}^{p_{i}}$ should be added or subtracted in order
to form the number $C$. Note that in order to have a possibility
to store $C$ in the computer memory, values $k$ and $m$ should be
finite.

Numbers $p_i$ in (\ref{3.13}) are  sorted in the decreasing order
with $ p_0=0$
\[
p_{m} >  p_{m-1}  > \ldots    > p_{1} > p_0 > p_{-1}  > \ldots
p_{-(k-1)}  >   p_{-k}.
 \]
They are called \textit{grosspowers} and they themselves can be
written in the form (\ref{3.13}).
 In the record (\ref{3.13}), we write
$\mbox{\ding{172}}^{p_{i}}$ explicitly because in the new numeral
positional system  the number   $i$ in general is not equal to the
grosspower $p_{i}$. This gives the possibility to write down
numerals without indicating grossdigits equal to zero.

The term having $p_0=0$ represents the finite part of $C$ because,
due to (\ref{3.2.1}), we have $c_0 \mbox{\ding{172}}^0=c_0$. The
terms having finite positive gross\-powers represent the simplest
infinite parts of $C$. Analogously, terms   having   negative finite
grosspowers represent the simplest infinitesimal parts of $C$. For
instance, the  number
$\mbox{\ding{172}}^{-1}=\frac{1}{\mbox{\ding{172}}}$ mentioned above
is infinitesimal. Note that all infinitesimals are not equal to
zero. Particularly, $\frac{1}{\mbox{\ding{172}}}>0$ because it is a
result of division of two positive numbers.

A number represented by a numeral in the form (\ref{3.13}) is called
\textit{purely finite} if it has neither infinite not infinitesimals
parts. For instance, 2 is purely finite and $2+3\G1^{-1}$ is not.
All grossdigits $c_i$ are supposed to be purely finite. Purely
finite numbers are used on traditional computers and for obvious
reasons have a special importance for applications.

All of the numbers introduced above can be grosspowers, as well,
giving thus a possibility to have various combinations of
quantities and to construct  terms having a more complex
structure. However, in this paper  we consider only purely finite
grosspowers. Let us give an example of multiplication of   two
infinite numbers $A$ and $B$ of this kind (for a comprehensive
description see \cite{Sergeyev_patent,informatica}).

\textbf{Example 2.} \label{e_m2}     Let us consider
  numbers $A$ and $B$, where
$$
A=14.3\mbox{\ding{172}}^{56.2} 5.4\mbox{\ding{172}}^{0},
\hspace{1cm}
B=6.23\mbox{\ding{172}}^{3}1.5\mbox{\ding{172}}^{-4.1}.
$$
The number $A$ has an infinite part  and a finite one. The number
$B$ has an infinite part and an infinitesimal one. Their product
$C$ is equal to
\[
\begin{tabular}{cr}\hspace {15mm}$C =  B \cdot A   =
 89.089\mbox{\ding{172}}^{59.2}21.45\mbox{\ding{172}}^{52.1}
33.642\mbox{\ding{172}}^{3}8.1\mbox{\ding{172}}^{-4.1}.$ &
\end{tabular} \hspace{10mm}
 \Box
 \]

We conclude this section by emphasizing that there exist different
mathematical languages and numeral systems and, if they have
different accuracies, it is not possible to use them together. For
instance, the usage of $`many$' from the language of Pirah\~{a}
in the record $4+ `many$' has no any sense because for Pirah\~{a} it
is not clear what is 4 and for people knowing what is 4 the accuracy
of the  answer `many' is too low. Analogously, the records of the
type $\G1 + \omega$, $\G1-\aleph_0 $, $\G1/\infty$, etc. have no
sense because they belong to languages developed   for different
purposes and having different accuracies.

\section{Numerical reconstruction of the   Taylor expansion of the solution on the Infinity Computer}

Let us return   to the problem (\ref{ode1}).  We suppose that   a
set of elementary functions ($a^{x}, \sin(x), \cos(x), $ etc.) is
represented at the Infinity Computer  by one of the usual ways
used in traditional computers (see, e.g. \cite{Muller}) involving
the argument $x$, finite constants, and four arithmetical
operations. Then the following theorem holds (the world
\textit{exact} in it means: with the accuracy of the computer
programme implementing $f(x,y)$ from (\ref{ode1})).

\begin{theorem}
\label{t_m1} Let us suppose that  for the solution $y(x),$ $ x \in
[a,b],$  of (\ref{ode1}) there exists  the  Taylor expansion
(unknown for us) and at  purely finite points $s \in [a,b],$ the
function $y(s)$ and all its derivatives assume purely finite
values   or are equal to zero. Then the Infinity Computer allows
one to reconstruct the Taylor expansion for $y(x)$ up to the
$k$-th derivative with exact values of $y'(x), y''(x), y^{(3)}(x),
\ldots y^{(k)}(x)$ after $k$ steps of the Euler method with the
step $h=\G1^{-1}$.
\end{theorem}

  \textbf{Proof.} Let us start  to execute on  the Infinite
Computer steps of the Euler method  following the rule
(\ref{ode2}) and using the infinitesimal step $h=\G1^{-1}$. Since
the problem (\ref{ode1}) has been stated using the traditional
finite mathematics, $x_0$ is purely finite. Without loss of
generality let us consider the first $k=4$ steps of the Euler
method (the value $k=4$ is sufficient to show the way of
reasoning; we shall use the formulae involved in this case later
in a numerical illustration). We obtain
 \beq
y_{1} = y_0 + \G1^{-1} f(x_0,y_0), \hspace{10mm} y_{2} = y_1 +
\G1^{-1} f(x_1,y_1),
 \label{ode9}
       \eeq
  \beq
y_{3} = y_2 + \G1^{-1} f(x_2,y_2),\hspace{10mm} y_{4} = y_3 +
\G1^{-1} f(x_3,y_3).
 \label{ode10}
       \eeq

The derivatives of the solution   $y(x)$ can be approximated in
different ways and with different orders of accuracy. Let us
consider approximations (see, e.g., \cite{Fornberg}) executed by
forward differences  $\triangle^j_{h}, 1 \le j \le k,$ with the
first order of accuracy and take $h=\G1^{-1}$  as follows
 \beq
  \triangle^k_{\tiny{\G1^{-1}}}=\sum^{k}_{i=0} (-1)^{i}
\left(\hspace{-1mm}
\begin{array}{c}
k \\
  i
  \end{array}\hspace{-1mm}   \right) y_{x_0+(k-i)\tiny{\G1^{-1}}}.
 \label{forward}
       \eeq
Then we have
  \beq
 y'(x_{0}) \approx \frac{ \triangle^1_{\tiny{\G1^{-1}}}  }{\G1^{-1} }
  + O\left(\G1^{-1} \right) =  \frac{y_{1} - y_{0}   }{\G1^{-1} }
  + O\left(\G1^{-1} \right),
\label{ode6.0}
       \eeq
 \beq
 y''(x_{0}) \approx \frac{ \triangle^2_{\tiny{\G1^{-1}}}   }{\G1^{-2} }
  + O\left(\G1^{-1} \right) =  \frac{y_{0} -2 y_{1} + y_{2}  }{\G1^{-2} }
  + O\left(\G1^{-1} \right),
\label{ode6}
       \eeq
 \beq
 y^{(3)}(x_{0}) \approx \frac{ \triangle^3_{\tiny{\G1^{-1}}}   }{\G1^{-3} }
  + O\left(\G1^{-1} \right) =  \frac{-y_{0} +3 y_{1} -3 y_{2} +
 y_{3}}{\G1^{-3} }
  + O\left(\G1^{-1} \right),
\label{ode6.1}
       \eeq
 \beq
 y^{(4)}(x_{0}) \approx \frac{ \triangle^4_{\tiny{\G1^{-1}}}   }{\G1^{-4} }
  + O\left(\G1^{-1} \right) =  \frac{y_{0} -4 y_{1} +
 6 y_{2} -4 y_{3} + y_{4}}{\G1^{-4} }
  + O\left(\G1^{-1} \right).
\label{ode6.2}
       \eeq

Since due to (\ref{ode1}) we can evaluate directly
$y'(x_{0})=f(x_0,y_0)$, let us start by considering the formula
(\ref{ode6}) (the cases with values of $k
> 2$ are studied by a complete analogy). Since $x_{0}$ is purely
finite, then due to our assumptions $y''(x_{0})$ is also purely
finite. This means that $y''(x_{0})$ does not contain infinitesimal
parts. Formula (\ref{ode6}) states that the error we have when
instead of $y''(x_{0})$ use its approximation
 \beq \widetilde{y''(x_{0})}  =
\frac{ \triangle^2_{\tiny{\G1^{-1}}}   }{\G1^{-2} }
 \label{ode8}
       \eeq
  is of the order
$\G1^{-1}$. The Infinity Computer works in such a way that it
collects different orders of \G1 in separate groups. Thus,
$\triangle^2_{\tiny{\G1^{-1}}}$ will be represented in the format
(\ref{3.13})
 \beq
 \triangle^2_{\tiny{\G1^{-1}}} =  c_{0}
\mbox{\ding{172}}^{0} +  c_{-1} \mbox{\ding{172}}^{-1} + c_{-2}
\mbox{\ding{172}}^{-2}  + \ldots  + c_{-m_2}
 \mbox{\ding{172}}^{-m_2}, \label{ode4}
       \eeq
where $m_2$ is a finite integer, its value depends on each concrete
$f(x,y)$ from (\ref{ode1}). Note that (\ref{ode4}) cannot contain
fractional grosspowers because the step $h=\G1^{-1}$ having the
integer grosspower $-1$ has been chosen in (\ref{ode9}),
(\ref{ode10}).

 It follows from (\ref{ode6}) and the
fact that $y''(x_{0})$ is purely finite that
$\widetilde{y''(x_{0})}$ contains a purely finite part and can
contain infinitesimal parts of the order $\G1^{-1}$ or higher.
This means that grossdigits $c_{0} = c_{-1}= 0$, otherwise after
division on $\G1^{-2}$ the estimate $\widetilde{y''(x_{0})}$ would
have infinite parts and this is impossible. Thus
$\widetilde{y''(x_{0})}$ has the following structure
 \beq
\widetilde{y''(x_{0})}  = c_{-2} \mbox{\ding{172}}^{0} +  c_{-3}
\mbox{\ding{172}}^{-1} + c_{-4} \mbox{\ding{172}}^{-2}  + \ldots +
c_{-m}
 \mbox{\ding{172}}^{-m+2}.
\label{ode7}
       \eeq
It follows from (\ref{ode6}) that $\widetilde{y''(x_{0})}$   can
contain an error  of the order $\G1^{-1}$ or higher. Since all the
derivatives of $y(x)$ are purely finite at $x_0$ and, in particular,
$y''(x_{0})$ is purely finite, the fact that  the finite part and
infinitesimal parts in (\ref{ode7}) are separated gives us that
$c_{-2}=y''(x_{0})$. Thus,   in order to have the exact value of
$y''(x_{0})$ it is sufficient to calculate
$\triangle^2_{\tiny{\G1^{-1}}}$ from (\ref{ode7}) and to take its
grossdigit $c_{-2}$ that will be equal to  $y''(x_{0})$.

 By a complete analogy the exact values of higher
derivatives can be obtained from (\ref{ode6.0}) -- (\ref{ode6.2})
and  analogous formulae using forward differences (\ref{forward}) to
approximate the $k$-th derivative $y^{(k)}(x_{0})$. It suffices just
to calculate $\triangle^k_{\tiny{\G1^{-1}}}$ and to take  the
grossdigit $c_{-k}$ that will be equal to the exact value of the
derivative $y^{(k)}(x_{0}). \hfill \Box $

Let us consider an illustrative numerical example. We emphasize that
the Infinity Computer solves it numerically, not symbolically, i.e.,
it is not necessary to translate the procedure implementing $f(x,y)$
in a symbolic form.

\textbf{Example 3.} \label{e_m3} Let us consider the problem
 \beq
 y'(x)=x-y, \hspace{1cm} y(0)=1,
 \hspace{1cm}
 \label{ode13}
       \eeq
taken from \cite{Adams}. Its exact solution is
  \beq
  y(x)=x-1+2e^{-x}.
 \label{ode11}
       \eeq
We start by applying formulae (\ref{ode9}) to calculate $y_{1}$ and
$y_{2}$:
 \[
y_{1} = 1 + \G1^{-1} \cdot (0-1) = 1 - \G1^{-1},
\]
\[
 y_{2} = 1 - \G1^{-1} +
\G1^{-1} (\G1^{-1}-1+\G1^{-1})= 1 - 2\G1^{-1} + 2\G1^{-2}.
 \]
 We have now the values $y_{0},\,\, y_{1},$ and $y_{2}$. Thus, we can
apply formula  (\ref{ode6}) and calculate
$\triangle^2_{\tiny{\G1^{-1}}}$ as follows
\[
  \triangle^2_{\tiny{\G1^{-1}}} =
 y_{0} -2 y_{1} + y_{2} = 1-2 + 2\G1^{-1} + 1 - 2\G1^{-1} +
 2\G1^{-2}=  2\G1^{-2}.
\]
Thus, $c_{-2}=2$. Let us now verify the obtained result and
calculate   the exact derivative $y''(0)$ using (\ref{ode11}). Then
we have $y''(x)=2e^{-x},  $ and $y''(0)=2$, i.e., $c_{-2}=y''(0)$.
Note that in this simple illustrative example $c_{-m}=0, \,\,  m>2,$
where $m$ is from (\ref{ode6}). In general, this is not the case and
$c_{-m}\neq 0$ can occur.

Let us proceed and calculate $y_{3}$ following (\ref{ode10}). We
have
  \[
y_{3} = 1 - 2\G1^{-1} + 2\G1^{-2} + \G1^{-1} (2\G1^{-1}-  1 +
2\G1^{-1} - 2\G1^{-2}  )= 1 - 3\G1^{-1} + 6\G1^{-2}-2\G1^{-3}.
\]
It then follows from (\ref{ode6.1}) that
 \[
\triangle^3_{\tiny{\G1^{-1}}} = -y_{0} +3 y_{1} -3 y_{2} +
 y_{3}  =
\]
 \[
  - 1 + 3 (1 - \G1^{-1}) -3 (1 - 2\G1^{-1} + 2\G1^{-2}) +
 1 - 3\G1^{-1} + 6\G1^{-2}-2\G1^{-3} =  -2\G1^{-3}.
\]
We can see that $c_{-3}=-2$. The exact derivative obtained from
(\ref{ode11}) is  $y^{(3)}(x)=-2e^{-x}$. As a consequence, we have
 $y^{(3)}(0)=-2$, i.e., $c_{-3}=y^{(3)}(0)$.

To calculate $y^{(4)}(0)$ we use (\ref{ode10}) and have
\[
 y_{4} = 1 - 3\G1^{-1} + 6\G1^{-2}-2\G1^{-3} +
\G1^{-1} (3\G1^{-1} -1 + 3\G1^{-1} - 6\G1^{-2}+2\G1^{-3})=
\]
\[
  1 - 4\G1^{-1} + 12\G1^{-2}- 8\G1^{-3} + 2\G1^{-4}.
\]
From (\ref{ode6.2}) we obtain
\[
\triangle^4_{\tiny{\G1^{-1}}} = y_{0} -4 y_{1} +
 6 y_{2} -4 y_{3} + y_{4} =  1 -4(1 - \G1^{-1})
 +6(1 - 2\G1^{-1} + 2\G1^{-2})-
\]
\[
 -4(1 - 3\G1^{-1} + 6\G1^{-2}-2\G1^{-3})
  +1 - 4\G1^{-1} + 12\G1^{-2}- 8\G1^{-3} + 2\G1^{-4}=
\]
\[
 -4(1 - 3\G1^{-1} + 6\G1^{-2}-2\G1^{-3})
  +1 - 4\G1^{-1} + 12\G1^{-2}- 8\G1^{-3} + 2\G1^{-4}=2\G1^{-4}.
\]
Again we obtain that $c_{-4}=y^{(4)}(0)=2$.

 Thus,   four
steps of the explicit Euler method with the infinitesimal step
$h=\G1^{-1}$ have been executed on the Infinity Computer. As a
result,   the first five exact items of the Taylor expansion of
$y(x)$ in the neighborhood of $x_0=0$ can be written:
 \beq
  y(x)=x-1+2e^{-x} \approx 1 - x + x^2 - \frac{x^3}{3} +
\frac{x^4}{12}.
 \label{ode12}
       \eeq
By a complete analogy it is possible to obtain additional terms in
the   expansion that correspond to higher derivatives of $y(x)$.

In Table~\ref{table1}, we present results  of experiments executed
(see \cite{Adams}) with the methods of Heun ($2^{d}$ order) and
Runge--Kutta ($4^{th}$ order) solving the problem (\ref{ode13}).
Both methods use $h=0.2$. Then we present results of the new methods
that execute first $k$ infinitesimal steps with $k$ going from 2 to
8 and then executing one finite step from the point $x_0=0$ to the
point $x=1$. The value $n_f$ is the number of evaluation of $f(x,y)$
executed by each method. The column $y_n$ shows the obtained
approximation at the point $x=1$ and the column $\varepsilon$ shows
the respective error $\varepsilon=y(1)-y_n$ where $n=5$ for the
methods of Heun   and Runge--Kutta and $n=1$ for the new methods.
\hfill
 $\Box$

\begin{table}[!t]
\caption{Comparison of   methods solving the problem (\ref{ode13})
where $n_f$ is the number of evaluation of $f(x,y)$ executed by a
method to obtain an approximated solution $y_n$ at the point $x=1$}
\begin{center}\scriptsize
\begin{tabular}{@{\extracolsep{\fill}}|c|c|c|c|}\hline
 Method  & $n_f$ & $y_n$ &    $\varepsilon$   \\
\hline
Heun, $h=0.2$ &$10$ &0.741480  & -0.005721         \\
Runge--Kutta ($4^{th}$ order), $h=0.2$ &$20$ &0.735770  & -0.0000116        \\
\hline
$y(x,0)=1 - x + x^2 $  &$2$ & 1  & -0.264241118       \\
$y(x,0)=1 - x + x^2 - \frac{x^3}{3}$   &$3$ & 0.6666666667  &\, 0.069092216     \\
$y(x,0)= 1 - x + x^2 - \frac{x^3}{3} + \frac{x^4}{12}$ &$4$ &0.75  & -0.014241118      \\
$y(x,0)= 1 - x + x^2 - \frac{x^3}{3} + \frac{x^4}{12}-\frac{x^5}{60}$ &$5$ &0.7333333333  & \, 0.002425549    \\
$y(x,0)= 1 - x + x^2 - \frac{x^3}{3} + \frac{x^4}{12}-\frac{x^5}{60}+\frac{x^6}{360}$ &$6$ &0.7361111111  & -0.000352229      \\
$y(x,0)= 1 - x + x^2 - \frac{x^3}{3} + \frac{x^4}{12}-\frac{x^5}{60}+\frac{x^6}{360}-\frac{x^7}{2520}$ &$7$ & 0.7357142857  &\, 0.000044597     \\
$y(x,0)= 1 - x + x^2 - \frac{x^3}{3} +
\frac{x^4}{12}-\frac{x^5}{60}+\frac{x^6}{360}-\frac{x^7}{2520}+\frac{x^8}{20160}$
&$8$ &0.7357638889  & -0.000005007
     \\
\hline
\end{tabular}
\end{center}
\label{table1}
\end{table}

In cases, where it is not possible to evaluate $f(x,y)$ at the
points $x_n+\G1^{-1}$, $x_n+2\G1^{-1}$,\, $x_n+3\G1^{-1}, \ldots$
(for instance, when we should solve the problem over an interval
$[a,b]$ and $x_n=b$) the following corollary can be useful.

\begin{corollary}\label{c_1}
Under conditions of Theorem~\ref{t_m1} the backward  differences
calculated at the points $x_n-\G1^{-1}$, $x_n-2\G1^{-1}$,\,
$x_n-3\G1^{-1}, \ldots ,  x_n-k\G1^{-1}$  can be used to calculate
the derivatives of $y(x)$ at the point $x_n$.
\end{corollary}

 \textbf{Proof.} The  backward difference (see e.g.,
\cite{Fornberg}) of the order $k$ with $h=\G1^{-1}$ is calculated as
follows
\[
  \nabla^k_{\tiny{\G1^{-1}}}=\sum^{k}_{i=0} (-1)^{i}
\left(\hspace{-1mm}
\begin{array}{c}
k \\
  i
  \end{array}\hspace{-1mm}   \right) y_{x_0-i\tiny{\G1^{-1}}}.
\]
The rest of the proof is completely analogous to the proof of
 the theorem and is so omitted. \hfill
 $\Box$

Thus, if the region of   interest $[a,b]$ from (\ref{ode1}) belongs
to the region of convergence of the Taylor expansion for the
solution $y(x)$ around the point $x_0$ then it is not necessary to
construct iterative procedures involving several steps with finite
values of $h$ and it becomes possible to calculate approximations of
the desired order by executing only one finite step.

\section{Examples  of the usage of infinitesimals in the new\\
computational framework}

The approach introduced in the previous section gives the
possibility to construct a variety of new numerical methods  for the
Infinity Computer by using both infinitesimal and finite values of
$h$. The general step $n$ of   a method of this kind for solving
(\ref{ode1}) can be described as follows:
 \bd
 \item (i) take the point $(x_n,y_n)$, choose
a value   $k_n$, and execute $k_n$ steps of the Euler method
starting from $x_n$ by using $h=\G1^{-1}$;
 \item (ii) calculate exact values of $y'(x),$ $
y''(x),$ $y^{(3)}(x), \ldots ,$ $y^{(k_n)}(x)$ at the point
$(x_n,y_n)$ following the rules described in Theorem~\ref{t_m1};
 \item (iii) construct the truncated Taylor
expansion of the order $k_n$;
 \item (iv) execute a single step from the point $x_n$ to $x_{n+1} = x_{n} +
h_i$ using  the constructed Taylor expansion and a finite value
of~$h_n$ (steps of the kind $h_n-\G1^{-1}$ or $h_n+\G1^{-1}$ can be
also used).
 \ed

The general step described above allows one to construct numerical
algorithms for solving (\ref{ode1}) by executing several iterations
of this kind.   Many numerical methods (see
\cite{Butcher,Henrici,Quarteroni}) can be used as a basis for such
developments. Due to the easy way  allowing us to  calculate exact
higher derivatives at the points $(x_n,y_n)$, methods that use
higher derivatives are of the main interest. The fact that to
increase the accuracy it is necessary just to execute one additional
infinitesimal step without performing additional finite steps (i.e.,
the whole work executed at a lower level of accuracy is used
entirely at a higher level of accuracy) is an additional advantage
and suggests to construct adaptive methods (for instance, if one
wishes to change the finite step from $h_1$ to $h_2 > h_1$ or $h_3 <
h_1$ then the same Taylor expansion can be used in all the cases). A
study of such methods will be done in a separate paper. Hereinafter,
since the usage of numerical infinitesimals is a new topic,   we
give a number of   examples  showing how the new computational
framework can be used in the context of numerical solving ODEs.

The rest of this section is organized as follows. In the first
subsection, we present two simple iterative methods using low
derivatives (a lower order of derivatives is used for expository
reasons). In the second subsection, we present a technique that can
be used to obtain an additional information with respect to
approximations of derivatives of the solution. In the last
subsection, we discuss how an automatic control of rounding errors
 can be executed during the evaluation of  $f(x,y)$ at the points
$(x_n,y_n)$.

\subsection{A simple method and possibilities of its improvements}

  We start by introducing the  \textit{Method 1} that uses only the
first and the second derivatives at each iteration to construct the
Taylor expansion by applying formulae (\ref{ode9}), (\ref{ode6.0}),
and (\ref{ode6}). Thus, this method at the current point $x_n$
executes twice the Euler step with $h=\G1^{-1}$ and then makes the
step with a finite value of $h$ by using the obtained Taylor
expansion. Therefore, during these three steps (two infinitesimal
steps and one finite step) the function $f(x,y)$ is evaluated twice
and only during the infinitesimals steps. Let us use the number step
$n$ to count the executed finite steps and denote by $y(x,z)$ the
Taylor expansion of the solution $y(x)$ calculated by the method
during the infinitesimal steps at the neighborhood of the point
$z=x_{n-1}$. Then, we can calculate $y_n$ as $y_n = y(h,x_{n-1})$
with a finite value of $h$.

\textbf{Example 4.} \label{e_m4} We test this method on the problem
(\ref{ode13}) with the finite step $h=0.2$ (see Table~\ref{table2}).
By applying the procedure described above with six digits after the
dot we have that
 \beq
  y(x,0) = 1 - x + x^2,   \hspace{1cm} y_1 = y(0.2,0) = 0.84,
  \label{ode14}
       \eeq
 \[
   y(x,0.2) = 0.84 - 0.64 x + 0.82 x^2, \hspace{1cm} y_2 = y(0.2,0.2) = 0.7448,
 \]
 \[
  y(x,0.4) = 0.7448 - 0.3448 x + 0.6724 x^2, \hspace{1cm} y_3 = y(0.2,0.4) =
  0.702736,
   \]
  \[
   y(x,0.6) = 0.702736 - 0.102736 x + 0.551368 x^2, \hspace{3mm} y_4 = y(0.2,0.6) =
  0.704244,
 \]
   \[
   y(x,0.8) = 0.704244 + 0.095756 x + 0.452122 x^2, \hspace{3mm} y_5 = y(0.2,0.8) =
  0.741480.
 \]

\begin{table}[!t] \caption{\small{Two versions of Method 1
 constructing the Taylor expansion using the values
of the first and second derivatives calculated during two
infinitesimal Euler steps at the points $(x_n,y_n)$}}
\begin{center}\scriptsize
\begin{tabular}{@{\extracolsep{\fill}}|c|c|c|c|c|c|c|}\hline
 &  & \multicolumn{2}{c|}{Method $1.0$} & \multicolumn{3}{c|}{Method $1.1$}  \\
\cline{3-7}$n$&$x_n$ &$y_n$ &$\varepsilon_n$   & $y^{c}_n$ & $c_n$ &$\varepsilon_n$  \\
\hline
0 &$0.0$ &1.000000     & \,\,0.000000  & 1.000000 & 0.000000 & \,0.000000    \\
1 &$0.2$ &0.840000  & -0.002538  &  0.839200 &  0.000800 & -0.001738    \\
\hline
2 &$0.4$ &0.744800  &-0.004160  & 0.743344 & 0.001456   &-0.002704 \\
3 &$0.6$ &0.702736  & -0.005113  & 0.700742   & 0.001994 &-0.003119 \\
\hline
4 &$0.8$ &0.704244  & -0.005586 &  0.701808   & 0.002436    & -0.003150 \\
5 &$1.0$ &0.741480  & -0.005721 &  0.738682   & 0.002798    & -0.002923 \\
\hline
\end{tabular}
\end{center}
\label{table2}
\end{table}

It can be seen from Table \ref{table2} that the results obtained by
the new method (see column Method~1.0) for the values $y_n$ coincide
(see \cite{Adams}) with the results obtained by applying the
modified Euler's method (called also Heun's method) that evaluates
$f(x,y)$ twice at each iteration as the new method does.

As it can be seen from the formulae above, the Method 1  at each
point $x_n$ provides us not only with the value $y_n$ but also with
the first and the second derivatives of $y(x)$ at the point
$(x_n,y_n)$. We shall denote them as $y'_n(x_{n})$ and
$y''_n(x_{n})$ where
\[
y'_n(x_{n})= y'(0,x_n),  \hspace{1cm} y''_n(x_{n})=y''(0,x_n).
\]
 For
instance, we have
\[
  y'(x,0.4) =   - 0.3448   + 1.3448 x, \hspace{1cm} y'_2(0.4)=-
  0.3448,
   \]
   \[
 \hspace{25mm}    y''(x,0.4) =
 1.3448, \hspace{1cm} y''_2(0.4)=1.3448. \hspace{20mm}
 \Box
   \]
Note that the values of derivatives calculated at the points
$(x_n,y_n)$ are exact in the sense that they provide values of the
derivatives for the solution with the initial condition $(x_n,y_n)$.
They can be used to estimate derivatives of $y(x)$ at the point
$(x_n,y(x_n))$.

 \begin{figure}[t]
  \begin{center}
    \epsfig{ figure = 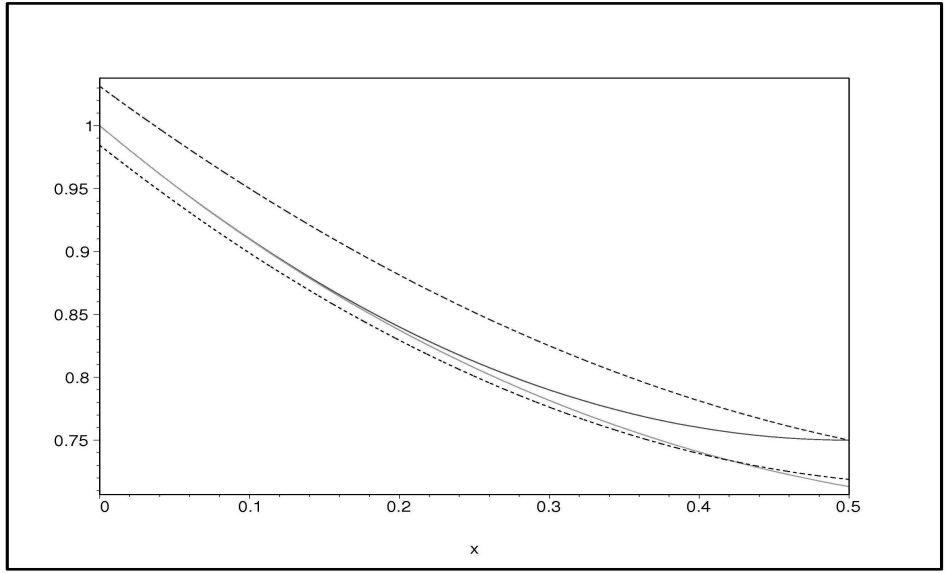, width = 4.9in, height = 3in, silent = yes }
    \caption{The graph of $y(x,0)$
is represented by the top (black) solid line and the graph of $y(x)$
is represented by the lower (grey) solid line (both $y(x,0)$ and
$y(x)$ start at the point $(0,1)$); the graph of the
 function $\bar{y}(x)$
is shown   by the top dashed line; the graph of the function
$r_1(x)$ is shown by the lower dashed line.}
 \label{figura_1}
  \end{center}
\end{figure}

The possibility to improve the accuracy of a solution going forward
and backward has been  the main idea of many numerical methods for
ODEs (see, e.g., \cite{Cash_Considine,Iavernaro_Mazzia}). Results
for the Method 1.0 presented in Table \ref{table2} can be improved
using similar ideas in the new framework using the values of
derivatives that the method provides during its functioning. Let us
consider the exact solution $y(x)$ from (\ref{ode11}) and its
approximation $y(x,0)$ from (\ref{ode14}) on the interval $[0,0.5]$
(we choose this interval that is larger then the step $h=0.2$ used
in Table~\ref{table2} in order to be able to show the idea more
clearly graphically). In Fig.~\ref{figura_1}, the graph of $y(x,0)$
is represented by the top (black) solid line and the graph of $y(x)$
is represented by the lower (grey) solid line. Naturally, both of
them start at the point $(0,1)$ and we then have that
\[
y(0.5)=0.713061319, \hspace{1cm} y(0.5,0)=0.75.
 \]
 Thus, the error $\varepsilon_1$ at the point $x_1=0.5$ is
  \beq
 \varepsilon_1 = y(0.5)- y(0.5,0) =
-0.036938681.
   \label{ode16.0}
       \eeq
By using   our equation (\ref{ode13}) and executing two
infinitesimal steps from the point $x_1=0.5$ to the points
$0.5+\G1^{-1}$ and $0.5+2\G1^{-1}$ we obtain
 \beq
 y(x,0.5) = 0.75-0.25  x + 0.625 x^2.
   \label{ode16}
       \eeq
In the Method 1.0, this formula would be used to go forward from the
point $x_1=0.5$ by executing a new finite step. Instead of this, let
us use this formula to go backward from the point $x_1=0.5$ to the
starting point $x_0=0$. Since in (\ref{ode16}) $x$ is intended to be
a removal from the point $x_1=0.5$, we need to take   into account
this fact. The graph of the obtained function
\[
 \bar{y}(x) = 0.75+0.25 (0.5-x)+ 0.625 (0.5-x)^2 = 1.03125-0.875x +0.625x^2
\]
is shown in Fig.~\ref{figura_1} by the top dashed line.

It can be seen that $y(x,0)$ from (\ref{ode14}) does not coincide
with the obtained function $\bar{y}(x)$. Let us construct a new
quadratic approximation $r_1(x)$ oriented on a better result to be
obtained at the point $x_1=0.5$. The function $r_1(x)$ is built
using coefficients from both $y(x,0)$ and $\bar{y}(x)$ by taking
their average with the weights $\frac{1}{2}$ and $\frac{1}{2}$ (the
general way to mix them is, obviously, $\tau$ and $1 - \tau,\, 0
<\tau<1$) as follows:
  \[
   r_1(x)= y(0,0)+ \frac{1}{2}(y(0,0)-\bar{y}(0))
   +\frac{1}{2}(y'(0,0)+\bar{y}'(0))x+\frac{1}{4}(y''(0,0)+\bar{y}''(0))x^2=
   \]
   \[
   1+\frac{1}{2}(1-1.03125)+\frac{1}{2}(-1-0.875)x+\frac{1}{2}(1+0.625)x^2=
   \]
    \beq
    0.984375-0.9375 x + 0.8125 x^2.
   \label{ode17}
       \eeq
In Fig.~\ref{figura_1}, the graph of the function $r_1(x)$ is shown
by the lower dashed line. The function $r_1(x)$ provides us the
value $r_1(0.5)=0.718750$ and the respective error
 \beq
 \overline{\varepsilon}_1 = y(0.5)- r_1(0.5) =
-0.005688681.
   \label{ode18}
       \eeq
that is better than the error $\varepsilon_1 = -0.036938681$ from
(\ref{ode16.0}) that is obtained by calculating $y(0.5,0)$.

The possibility to calculate  corrections to approximations opens
the doors to various modifications. For instance, it is possible to
execute two additional infinitesimal steps at the point $x_1=0.5$
using the value $r_1(0.5)$ instead of $y(0.5,0)$. In general, this
means that instead of setting $y_n=y(x_{n},x_{n-1})$ as it is done
by the Method 1.0 we put $y_n=r_n(x_n)$. Obviously, this means
 that it is necessary to evaluate $f(x,y)$ two times more with
 respect to the Method 1.0. Otherwise it is also possible to use the corrected
 value $r_n(x_n)$ with the derivatives that have been calculated for
 $y(x_{n},x_{n-1})$.

 Another possibility would be   the use of the functions $y(x,x_n)$ at each
 point~$x_n$, i.e., to put $y_n=y(x_{n},x_{n-1})$,  and to calculate the
global correction following the rule
 \beq
c_n=c(x_n)= c(x_{n-1}) + r_n(x_n)-y(x_{n},x_{n-1}),
  \label{ode19}
       \eeq
       starting from the first correction (in our example
       $c(x_1)=c(0.5)=
       0.031250$)
   \[
c(x_1)=   r_1(x_1)-y(x_{1},x_{0}).
\]
In this way we can approximate the exact solution $y(x_n)$ by the
corrected value
 \[
y^{c}_n = y(x_{n},x_{n-1})+c(x_n).
 \]
In Table~\ref{table2}, results for this algorithm are presented in
the column \emph{Method~1.1} where the error $\varepsilon_n$ is
calculated as
\[
\varepsilon_n=y(x_n)-y^{c}_n.
\]
 Notice   that the
correction obtained at the final point has been calculated using
Corollary~1.

We conclude this subsection by a reminder that Theorem~\ref{t_m1}
gives us the possibility to easily construct   higher-order methods.
Two methods described above just show examples of the usage of
infinitesimals for building algorithms for solving ODEs.

\subsection{Approximating derivatives of the solution}

In this subsection, we show how approximations of derivatives at the
point $x_n$ can be obtained using the information calculated at the
point $x_{n-1}$. For this purpose, instead of the usage of a finite
step $h$, the steps $h-\G1^{-1}$ or $h+\G1^{-1}$ can be used. To
introduce this technique we need to recall the following theorem
from \cite{Num_dif}.

\begin{theorem}
\label{t_m2} Suppose that: (i) for a function $s(x)$   calculated by
a procedure implemented at the Infinity Computer  there exists an
unknown  Taylor expansion in a finite neighborhood $\delta(z)$ of a
purely finite point $z$; (ii) $s(x),$ $s'(x), s''(x), \ldots
s^{(k)}(x)$ assume purely finite values or are equal to zero at
purely finite $x \in \delta(z)$; (iii) $s(x)$ has been evaluated at
a point $z+\mbox{\ding{172}}^{-1} \in \delta(z)$. Then the Infinity
Computer returns the result of this evaluation in the positional
numeral system with the infinite radix~\ding{172} in the following
form
 \beq
  s(z+\mbox{\ding{172}}^{-1}) = c_{0}
\mbox{\ding{172}}^{0}   c_{-1} \mbox{\ding{172}}^{-1}  c_{-2}
\mbox{\ding{172}}^{-2}   \ldots c_{-(k-1)}
 \mbox{\ding{172}}^{-(k-1)} c_{-k}
 \mbox{\ding{172}}^{-k},
\label{m1}
       \eeq
      where
 \beq
   s(z) = c_{0}, \,\,
s'(z) =    c_{-1}, \,\, s''(z)= 2!   c_{-2}, \,\,
 \ldots \,\,
 s^{(k)}(z)=k!  c_{-k}.
       \label{m2}
       \eeq
\end{theorem}

The theorem tells us that if we take a purely finite point $z$ and
evaluate on the Infinity Computer $s(x)$ at the point $z+\G1^{-1}$
then from the computed $s(z+\G1^{-1})$ we can easily extract $s(z)$,
$s'(z)$, $s''(z)$, etc. To apply this theorem to our situation  we
can take as $s(x)$ the Taylor expansion for $y(x)$ constructed up to
the $k$th derivative using infinitesimal steps $\G1^{-1}$,
$2\G1^{-1}$, $ \ldots, k\G1^{-1}$. Then, if we take as $z$ a purely
finite step $h$ and evaluate on the Infinity Computer $s(x)$ at the
point $h+\G1^{-1}$ then we obtain $s(h)$, $s'(h)$, $s''(h)$, etc.

For instance, let us take $s(x)=s_2(x)=y(x,0)$, where $y(x,0)$ is
from (\ref{ode14}) and $s_2(x)$ indicates that we use two
derivatives in the Taylor expansion. Then, we have
 \[
 s_2(0.2+\G1^{-1}) = 1 - (0.2+\G1^{-1})+ (0.2+\G1^{-1})^2 =
 0.84 - 0.6\G1^{-1}+\G1^{-2}.
 \]
We have the exact (where the word ``exact'' again means: with the
accuracy of the implementation of $s(x)$) values $s(0.2)=0.84$,
$s'(0.2)=-0.6$, $s''(0.2)=1$ for the function $s(x)$. These values
can be used to approximate the respective values $y(0.2)$,
$y'(0.2)$, $y''(0.2)$ we are interested in. Moreover, we can
adaptively obtain an information on the accuracy of our
approximations by consecutive improvements. If we calculate now
$y^{(3)}(0)$ from (\ref{ode6.1}) then we can improve our
approximation by setting $s(x)=s_3(x)$ where
\[
s_3(0.2+\G1^{-1})= s_2(0.2+\G1^{-1}) - \frac{1}{3}(0.2+\G1^{-1})^3 =
 \]
 \[
s_2(0.2+\G1^{-1}) - 0.002667 - 0.04 \G1^{-1}-0.2\G1^{-2}
-\frac{1}{3}\G1^{-3} =
\]
 \beq
    0.837333 - 0.64 \G1^{-1}+0.8\G1^{-2}
-\frac{1}{3}\G1^{-3}.
  \label{ode15}
       \eeq
  Note, that to obtain this information
we have calculated only the additional part  of $s_3(0.2)$ taking
the rest from the already calculated value $s_2(0.2)$.

 Analogously, if we calculate now
$y^{(4)}(0)$ from (\ref{ode6.2}) then we can improve our
approximation again by setting $s(x)=s_4(x)$ where
\[
s_4(0.2+\G1^{-1})= s_3(0.2+\G1^{-1}) + \frac{1}{12}(0.2+\G1^{-1})^4
=
 \]
 \[
 s_3(0.2+\G1^{-1}) +  0.000133 +
0.002667 \G1^{-1}+0.02\G1^{-2}
+0.066667\G1^{-3}+\frac{1}{12}\G1^{-4} =
\]
 \[
    0.837466 - 0.637333 \G1^{-1}+0.82 \G1^{-2}
-0.266667\G1^{-3}+\frac{1}{12}\G1^{-4}.
  \]
 Since we
have used the convergent Taylor expansion of the fourth order, the
errors in calculating $y(0.2)$, $y'(0.2)$, and $y''(0.2)$ are of the
orders 5, 4, and 3, respectively.

\subsection{An automatic control of rounding errors}

In the previous sections, we have supposed that the evaluation of
derivatives of $y(x)$ was done exactly, i.e., the procedure for
evaluating $f(x,y)$ was sufficiently precise and it was possible to
neglect rounding errors. The executed numerical examples presented
above   satisfied this assumption. Let us now study what the
Infinity Computer can give us when $f(x,y)$ is calculated with
errors. Hereinafter we suppose again that for the solution $y(x),$ $
x \in [a,b],$ of (\ref{ode1}) there exists  the  Taylor expansion
(unknown for us) and at  purely finite points $s \in [a,b],$ the
function $y(s)$ and all its derivatives assume purely finite values
or are equal to zero. In addition, we assume that the same
conditions hold for all approximations of $y(x)$ the method will
deal with.

Let us consider   formulae (\ref{ode9}) for calculating $y_{1}$ and
$y_{2}$ together with formulae (\ref{ode6.0}) and (\ref{ode6}) used
for approximating    $ y'(x_{0})$ and $ y''(x_{0})$. Suppose that
$f(x_0,y_0)$ is calculated with an unknown error $\epsilon_1$.  Then
we have that instead of the derivative $y'(x_{0})$ and the point
$y_1$  we have
 \beq
\tilde{y}'(x_{0})= f(x_0,y_0)- \epsilon_1 = y'(x_{0})- \epsilon_1
 \label{ode21}
       \eeq
       \[
\tilde{y}_{1} = y_0 + \G1^{-1} (f(x_0,y_0)- \epsilon_1)= y_1 -
\epsilon_1\G1^{-1},
\]
Analogously, calculation of the point $y_{2}$ will give us
$\tilde{y}_{2}$ and $\triangle^2_{\tiny{\G1^{-1}}}$ will be
calculated with errors, as well, giving us $
\tilde{\triangle}^2_{\tiny{\G1^{-1}}}$. Let us study the structure
of this forward difference
 \beq
  \tilde{\triangle}^2_{\tiny{\G1^{-1}}}= y_{0} -2 \tilde{y}_{1} + \tilde{y}_{2}= \tilde{y}_{2}-
  \tilde{y}_{1} -( \tilde{y}_{1} -y_{0})= \tilde{y}_{2}-
  \tilde{y}_{1} -( y_{1}
 -y_{0}) +  \epsilon_1 \G1^{-1} .
 \label{ode22}
       \eeq
By applying the argumentation analogous to that  used in
Theorem~\ref{t_m1} together with our assumptions on purely
finiteness of all derivatives of approximations of $y(x)$ we have
that
 \beq
\tilde{y}_{2}-
  \tilde{y}_{1} -( y_{1}
 -y_{0}) =  \widetilde{c}_{-2}
\mbox{\ding{172}}^{-2}  + \ldots  + \widetilde{c}_{-m_2}
 \mbox{\ding{172}}^{-m_2},
 \label{ode23}
       \eeq
where the coefficients $\widetilde{c}_{-2}, \ldots  ,
\widetilde{c}_{-m_2}$ are affected by rounding errors and errors
incorporated in $\tilde{y}_{1}$ and $\tilde{y}_{2}$. Thus, instead
of the exact second derivative that has been obtained in
Theorem~\ref{t_m1} from the coefficient $c_{-2}$ of $\G1^{-2}$, the
coefficient $\widetilde{c}_{-2}$ gives us an approximation
$\tilde{y}''(x_{0})$ of $  y''  (x_{0})$, namely,
\[
\widetilde{c}_{-2}=\tilde{y}''(x_{0}) = y''(x_{0})-\epsilon_2,
\]
where $\epsilon_2$ is an error we have got during the calculation of
$y''(x_{0})$.

 Let us rewrite now (\ref{ode22}) in the decreasing
orders of the powers of grossone using the representation
(\ref{ode23}), i.e., as the Infinity Computer does it. We have
 \beq
  \tilde{\triangle}^2_{\tiny{\G1^{-1}}}= \epsilon_1 \G1^{-1}  +   \widetilde{c}_{-2}
\mbox{\ding{172}}^{-2}  + \ldots  + \widetilde{c}_{-m_2}
 \mbox{\ding{172}}^{-m_2}.
 \label{ode24}
       \eeq
This means that by calculating
$\tilde{\triangle}^2_{\tiny{\G1^{-1}}}$ we have obtained also the
error $\epsilon_1$ that we have got at the previous infinitesimal
step (see (\ref{ode21})). We are able now to reestablish the exact
value of the first derivative $y'(x_{0})$ using the approximative
value $\tilde{y}'(x_{0})$ calculated in (\ref{ode21}) and the
grossdigit corresponding to $\G1^{-1}$ by taking it from
$\tilde{\triangle}^2_{\tiny{\G1^{-1}}}$ in (\ref{ode24}), i.e., we
have
 \[
y'(x_{0}) = \tilde{y}'(x_{0})+ \epsilon_1.
 \]
By a complete analogy we can continue and calculate
 \beq
\tilde{\triangle}^3_{\tiny{\G1^{-1}}}= -\epsilon_1 \G1^{-1}  +
\epsilon_2 \G1^{-2}  +  \widetilde{c}_{-3} \mbox{\ding{172}}^{-3}  +
\ldots  + \widetilde{c}_{-m_3}
 \mbox{\ding{172}}^{-m_3},
   \label{ode26}
       \eeq
 \[
y''(x_{0}) = \tilde{y}''(x_{0})+ \epsilon_2.
 \]
Note that in (\ref{ode26})  $\epsilon_1$ (that can be either
positive or negative) appears with the alternated sign following the
formulae of forward differences. In fact, in
$\tilde{\triangle}^3_{\tiny{\G1^{-1}}}$ we have $ y_{1}
 -y_{0}$ whereas in
$\tilde{\triangle}^2_{\tiny{\G1^{-1}}}$ we have $-( y_{1}
 -y_{0})$. Analogously, the same alternation happens for higher
derivatives.

 In general, in order to
calculate the $(k-1) {th}$ derivative $y^{(k-1)}(x_{0})$ it is
necessary to calculate the approximation $\tilde{y}^{(k-1)}(x_{0})=
\widetilde{c}_{-(k-1)}$ and then to extract the error
$\epsilon_{k-1}$ (that can be negative or positive ) from
 \[
\tilde{\triangle}^k_{\tiny{\G1^{-1}}}= (-1)^{k}\epsilon_1 \G1^{-1} +
\ldots (-1)^{k-i-1}\epsilon_i \G1^{-i} + \ldots
 \]
 \[ - \epsilon_{k-2}
\G1^{-(k-2)} + \epsilon_{k-1} \G1^{-(k-1)} +
 \widetilde{c}_{-k}
\mbox{\ding{172}}^{-k} + \ldots + \widetilde{c}_{-m_k}
 \mbox{\ding{172}}^{-m_k},
 \]
 \beq
y^{(k-1)}(x_{0}) = \tilde{y}^{(k-1)}(x_{0})+ \epsilon_{k-1}.
  \label{ode25}
       \eeq
If there exists an index $j, \,\, 1 \le j < k,$ such that
$\epsilon_{1}=\ldots\epsilon_{j}=0$, then $y^{(k-1)}(x_{0})$ is
calculated again by the formula (\ref{ode25}) but it follows
 \[
\tilde{\triangle}^k_{\tiny{\G1^{-1}}}= (-1)^{k-j-1}\epsilon_{j+1}
\G1^{-(j+1)} + \ldots
 \]
 \beq
 + \epsilon_{k-1}
\G1^{-(k-1)}  + \widetilde{c}_{-k} \mbox{\ding{172}}^{-k} + \ldots
+ \widetilde{c}_{-m_k}
 \mbox{\ding{172}}^{-m_k}.
 \label{ode28}
       \eeq

Thus, either $f(x,y)$ is evaluated exactly or rounding errors are
present, the Infinity Computer is able to calculate the derivatives
of the solution exactly. Let us illustrate the theoretical results
presented above by a numerical example.

\textbf{Example 5.} \label{e_m5} Let us consider the following test
problem\footnote{The author thanks Prof. H. P. Langtangen for
drawing the author's attention to this nice example.} taken from
\cite{Langtangen}. The ODE
 \beq
y'(x) = -\frac{x-c}{s^2} (y-1)
 \label{ode20}
       \eeq
 has the  exact solution that is the following Gaussian function
 \beq
 u(x) = 1 + e^{-\frac{1}{2}\left(\frac{x-c}{s}\right)^2}
 \label{ode27}
       \eeq
centered around $t=c$ and with characteristic width (standard
deviation) $s$. The initial condition is taken as the exact value
$y(0)=u(0)$ and the parameters are taken as $c=3, s=0.5$.
\begin{table}[!t]
\caption{Calculating approximations for derivatives $y^{(i)}(0), 1
\le i \le 12,$ for the problem (\ref{ode20}) at the point $(0,y(0))$
by applying the automatic control of rounding errors $\varepsilon_i$
and providing so the final accuracy $\delta_i$ }
\begin{center}\scriptsize
\begin{tabular}{@{\extracolsep{\fill}}|c|c|c|c|c|}\hline
$i$ & $\tilde{y}^{(i)}(0)$ & $\varepsilon_i$ &  $y^{(i)}(0)$   & $\delta_i$    \\
\hline
1 & $0.182759757\cdot 10^{-6}$ &  \hspace{-8mm}0.0000000000                & $0.182759757\cdot 10^{-6}$  & $-0.60449198\cdot 10^{-28} $         \\
2 & $0.207127725\cdot 10^{-5}$ & $0.609199190\cdot 10^{-7}$  & $0.213219716\cdot 10^{-5}$  & $-0.69190731\cdot 10^{-27}$    \\
3 & $0.233932489\cdot 10^{-4}$ & $0.731039028\cdot 10^{-6}$ &  $0.241242879\cdot 10^{-4}$  & $-0.78192941\cdot 10^{-26}$    \\
4 & $0.254888941\cdot 10^{-3}$ & $0.901614801\cdot 10^{-5}$  & $0.263905089\cdot 10^{-3}$  & $-0.83928642\cdot 10^{-25}$    \\
5 & $0.266975453\cdot 10^{-2}$ & $0.111117932\cdot 10^{-3}$  & $0.278087246\cdot 10^{-2}$  & $-0.85899500\cdot 10^{-24}$    \\
6 & $0.267228880\cdot 10^{-1}$ & $0.136947978\cdot 10^{-2}$  & $0.280923676\cdot 10^{-1}$  & $-0.82216871\cdot 10^{-23}$    \\
7 & \hspace{-2mm}$0.253489245\cdot 10^{0}$  & $0.168782291\cdot 10^{-1}$  & \hspace{-2mm}$0.270367474\cdot 10^{0}$   & $-0.71132365\cdot 10^{-22}$    \\
8 & \hspace{-2mm}$0.224980672\cdot 10^{1}$  & \hspace{-2mm}$0.208016667\cdot 10^{0}$  & \hspace{-2mm}$0.245782339\cdot 10^{1}$   & $-0.50790463\cdot 10^{-21}$    \\
9 & \hspace{-2mm}$0.182784086\cdot 10^{2}$  & \hspace{-2mm}$0.256371293\cdot 10^{1}$  & \hspace{-2mm}$0.208421215\cdot 10^{2}$   & $-0.19195037\cdot 10^{-20}$    \\
10 & \hspace{-2mm}$0.13002719\cdot 10^{3}$  & \hspace{-2mm}$0.315966218\cdot 10^{2}$  & \hspace{-2mm}$0.161623816\cdot 10^{3}$   & $ \hspace{2mm}0.25832403\cdot 10^{-19}$    \\
11 & \hspace{-2mm}$0.71638662\cdot 10^{3}$  & \hspace{-2mm}$0.389414314\cdot 10^{3}$  & \hspace{-2mm}$0.110580093\cdot 10^{4}$   & $ \hspace{2mm}0.86669850\cdot 10^{-18}$    \\
12 & \hspace{-2mm}$0.13588050\cdot 10^{4}$  & \hspace{-2mm}$0.479935826\cdot 10^{4}$  & \hspace{-2mm}$0.615816329\cdot 10^{4}$   & $ \hspace{2mm}0.16535675\cdot 10^{-16}$    \\
 \hline
\end{tabular}
\end{center}
\label{table3}
\end{table}

In Table \ref{table3}, we calculate derivatives $y^{(i)}(0), 1 \le i
\le 12,$ using (\ref{ode25})  with 30 digits in the mantissa, in
order to be able to catch the final accuracy $\delta_k$ presented in
the last column. It shows the final error obtained by subtracting
from the derivatives calculated using the explicit solution
(\ref{ode27}) and the derivatives $y^{(i)}(0), 1 \le i \le 12,$
i.e.,
 \[
 \delta_i = u^{(i)}(0)- y^{(i)}(0), \hspace{8mm} 1 \le i \le 12.
 \]
Let us make a few remarks regarding Table~\ref{table3}. First, at
$x=0$   with $c=3$ and $s=0.5$ it follows from (\ref{ode20}) that
$-\frac{x-c}{s^2} = 12$. In order to illustrate the situation
(\ref{ode28}), we have calculated $y_1$ using (instead of the
original expression $12(y_0-1)$ from (\ref{ode20}) leading to
$\varepsilon_1 \neq 0$) the expression $12 y_0-12$  that provides
$\varepsilon_1 = 0$ when it is used in $ y_1= y_0 + \G1^{-1}(12
y_0-12)$.

Then, it is worthwhile to notice that almost through the whole
Table~\ref{table3} (and in spite of large values of higher
derivatives)   the relative error has the constant order equal to
$10^{-22}$ (it can be easily seen from Table~\ref{table3} that
$m-n=22$ where $m$ is the exponent of $y^{(i)}$ and $n$ is the
exponent of $\delta_i$). Notice also that at the last line of the
Table the error $\varepsilon_{12}$ is even larger than the
approximation $\tilde{y}^{(12)}(0)$.
  \hfill $\Box$

\section{Conclusion}

In this paper, a new framework for solving ODEs has been introduced.
The new approach allows us to work numerically not only with usual
finite numbers but also with different infinitesimal and infinite
values on a new kind of a computational device called the Infinity
Computer (it has been patented and its working prototype exists).
The structure of numbers we work on the new computer is   more
complex and, as a result, we face new computational possibilities.
In particular, the presence of different numerical infinitesimals
makes it possible to use infinitesimal steps for solving ODEs. The
following results have been established in the new framework.

i. It has been shown that (under the assumption that the person
solving the ODE does not know the structure of $f(x,y)$, i.e., it is
a ``black box'' for him/her) the Infinity Computer is able   to
calculate numerical values of the derivatives of $y(x)$ of the
desired order without the necessity of  an analytical (or
symbolical) computation of the respective derivatives by the
successive derivation of the ODE as it is usually done when the
Taylor method is applied.

ii. If the region of our interest $[a,b]$ belongs to the region of
convergence of the Taylor expansion for the solution $y(x)$ in the
neighborhood of the point $x_0$, then it is not necessary to
construct iterative procedures involving several steps with finite
values of $h$. It becomes possible to calculate approximations of
the desired order $k$ by executing $k$ infinitesimal steps and only
one finite step.

iii. Approximations of derivatives of $y(x)$ at the point $x_n$ can
be obtained using the information calculated at the point $x_{n-1}$.
For this purpose, instead of the usage of a finite step $h$, the
steps $h-\G1^{-1}$ or $h+\G1^{-1}$ can be used. Methods going
forward and backward and working with approximations of derivatives
can be proposed in the new framework.

iv. The last subsection of the manuscript shows that either $f(x,y)$
is evaluated exactly or rounding errors are present, the Infinity
Computer is able to perform,    by means of a smart usage of
infinitesimals, an automatic control of the accuracy of computation
of the derivatives of the solution.

v. Theoretical results have been illustrated by a number of
numerical examples.

\end{document}